\documentclass[11pt]{article}

      \usepackage{latexsym}
      \usepackage[centertags]{amsmath}
\usepackage{amsfonts}
\usepackage{amssymb}

\makeatletter\@addtoreset{equation}{section} \makeatother

\renewcommand\thefigure{\thesection.\@arabic\c@figure}
\renewcommand\thetable{\thesection.\@arabic\c@table}

      \newtheorem{theorem}{Theorem}[section]

      \newtheorem{lemma}[theorem]{Lemma}
      \newtheorem{Lem}{Lemma}[section]

      \newtheorem{Prop}[theorem]{Proposition}

      \def\nn{\nonumber}
      \def\rf#1{\mbox{$(\ref{#1})$}}

      \def\be{\begin{equation}} 
      \def\ee{\end{equation}} 
      \def\beqn{\begin{eqnarray}} 
      \def\eeqn{\end{eqnarray}} 
      \def\beq{\begin{eqnarray*}} 
      \def\eeq{\end{eqnarray*}}
      \def\proof{{\noindent\bf Proof\quad}\ }
      \def\mb{\mbox} 
      \def\ga{\gamma} 
      \def\la{\lambda} 
      \def\ra{\rightarrow} 

\begin{document}
      \title{Asymptotic Results for the Two-parameter Poisson-Dirichlet Distribution }

 \author{Shui Feng\thanks{Research supported by
      the Natural Science and Engineering Research Council of Canada}\\Department of Mathematics \\and Statistics\\McMaster
      University\\Hamilton, Ontario\\
      Canada L8S 4K1\\
      shuifeng@mcmaster.ca\\ \and Fuqing Gao\thanks{Research supported  by the NSF of
China(No.10571139).}\\School of Mathematics and Statistics\\Wuhan
University\\ Wuhan 430072, China\\fqgao@whu.edu.cn}

\date{}
\maketitle

\begin{abstract}
The two-parameter Poisson-Dirichlet distribution is the law of a
sequence of decreasing nonnegative random variables with total sum
one. It can be constructed from stable and Gamma subordinators with
the two-parameters, $\alpha$ and $\theta$, corresponding to the
stable component and Gamma component respectively. The moderate
deviation principles are established for the two-parameter
Poisson-Dirichlet distribution and the corresponding homozygosity
when $\theta$ approaches infinity, and the large deviation principle
is established for the two-parameter Poisson-Dirichlet distribution when both $\alpha$ and $\theta$ approach zero.

\end{abstract}
\vspace*{.125in} \noindent {\bf Key words:}  Poisson-Dirichlet
      distribution, two-parameter Poisson-Dirichlet distribution, GEM representation,
      homozygosity, large deviations, moderate deviations.
      \vspace*{.125in}

      \noindent {\bf AMS 1991 subject classifications:}
      Primary: 60F10; Secondary: 92D10.

\section{Introduction}
\setcounter{equation}{0}

For $\alpha$ in $(0,1)$ and $\theta >
-\alpha$, let $U_k,
      k=1,2,\cdots$, be a sequence of independent
      random variables such that $U_k$ has $Beta(1-\alpha,\theta+ k\alpha)$ distribution.
      Set
      \be \label{GEM1}
      X^{\alpha,\theta}_1 = U_1,\  X^{\alpha,\theta}_n = (1-U_1)\cdots (1-U_{n-1})U_n,\  n \geq 2.
      \ee

      Then with probability one $$\sum_{k=1}^{\infty}X^{\alpha,\theta}_k =1,$$ and the law of
      $(X^{\alpha,\theta}_1,X^{\alpha,\theta}_2,\cdots)$
      is called the two-parameter GEM distribution.

      Let ${\bf P}(\alpha,\theta)=(P_1(\alpha,\theta), P_2(\alpha,\theta),\cdots)$ denote
      the descending order statistic of $(X^{\alpha,\theta}_1,X^{\alpha,\theta}_2,\cdots)$.
      The law of ${\bf P}(\alpha,\theta)$ is called the two-parameter
      Poisson-Dirichlet distribution and is denoted by
      $\Pi_{\alpha,\theta}$. The well-known one-parameter Poisson-Dirichlet
      distribution corresponds to $\alpha=0$.

      For each integer $m \geq 2$, taking a random sample of size $m$ from
a population with the two-parameter Poisson-Dirichlet distribution.
Given the population proportion ${\bf p}=(p_1,p_2,\ldots)$, the
probability that all samples are of the same type is given by
\[
H_m({\bf p})=\sum_{i=1}^{\infty}p_i^m,
\]
which is referred to as the homozygosity of order $m$.

      The main properties of the two-parameter Poisson-Dirichlet distribution are studied in
      Pitman and Yor \cite{PitmanYor97} including relations to subordinators, Markov chains, Brownian motion and
      Brownian bridges.  The detailed calculations of moments and parameter estimations were carried out in
      Carlton \cite{Car99}. In \cite{Der97} and the references therein one can find
      connections between
      two-parameter Poisson-Dirichlet distribution and models in physics including mean-field spin glasses, random map
      models,
      fragmentation, and returns of a random walk to origin. The two-parameter Poisson-Dirichlet distribution
      has also been used in
      macroeconomics and finance (\cite{Aoki06}).

      Many properties of the one-parameter Poisson-Dirichlet distribution have generalizations
      in the two-parameter setting including but not limited to the sampling
      formula (cf. \cite{Ewen72}, \cite{pitman92}), the
      Markov-Krein identity (cf. \cite{diaco96}, \cite{tsi97}), and
      subordinator representation (cf. \cite{Kingman75},
      \cite{PitmanYor97}). Recently, a large deviation principle (henceforth, LDP)is
      established in \cite{Feng07} for the two-parameter
      Poisson-Dirichlet distribution when $\theta$ goes to infinity.
      This is a generalization to the LDP result for the one-parameter Poisson-Dirichlet distribution in
      \cite{dawson-feng06}. Our first result here establishes the corresponding moderate
      deviation principle (henceforth, MDP). This can be viewed as a
      generalization of the MDP result in \cite{FengGao08} to the
      two-parameter setting. The MDP for the
      homozygosity is also established generalizing corresponding
      result in \cite{FengGao08}. In order to apply the Campbell's theorem, we turn to a representation of
      the two-parameter Poisson-Dirichlet distribution obtained in
      \cite{pitman92}.

      When $\alpha =0$, the one-parameter Poisson-Dirichlet distribution converges to $\delta_{(1,0,\ldots)}$
      as $\theta $ goes to zero.
      The corresponding LDP is established in \cite{Feng08} where a
      structure called ``energy ladder" is revealed. Our second main
      result generalizes this result to the two-parameter
      Poisson-Dirichlet distribution when both $\alpha$ and $\theta$
      go to zero. It turns out that the large deviation speed
      will depend on $\alpha$ if it converges to zero at a slower
      speed than that of $\theta$.

  The current paper is organized as follows.  Distributional results are derived in Section 2 using the
  change of measure formula and the subordinator representation. Section 3 is dedicated to establishing the MDP for
  $\Pi_{\alpha,\theta}$ when $\theta$ goes to infinity. The large $\theta$ MDP for the homozygosity is established
  in Section 4.
  In Section 5 we prove the LDP for $\Pi_{\alpha,\theta}$ when both $\alpha$ and $\theta$ go
  to zero.

\section{Marginal Distributions}

In this section, we derive the marginal distributions of the
two-parameter Poisson-Dirichlet distribution. The basic tools are
the change of measure formula and the subordinator representation.
For general concepts and theorems on MDP and LDP, we will refer to
\cite{dembo-zeitiouni98}.

>From now on, the parameter $\theta$ will be assumed to be positive
and $\alpha$ is in $(0,1)$. Let $\{\rho_s, s\geq 0\}$ be a
subordinator, an increasing process with stationary
independent increment, with no drift component. The Laplace
transform of $\rho_s$ is then given by

 \be E\left(\exp(-\la
\rho_s)\right)=\exp\left\{s\int_0^{\infty}(e^{-\la
x}-1)\Lambda(dx)\right\}, \la \geq 0, \ee where $\Lambda$ is the
L\'evy measure on $(0,+\infty)$ describing the distribution of the
jump sizes. Let $V_1(\rho_s)\geq V_2(\rho_s)\geq \cdots$ denote the
jump sizes of $\{\rho_s,s \geq 0\}$ over $[0,s)$ in decreasing
order.

If
\[
\Lambda(dx)=c_{\alpha} x^{-(1+\alpha)}dx.
\]
for some $c_{\alpha} >0$, then the subordinator is called a stable subordinator with index $\alpha$ and is denoted by
$\{\tau_s,s\geq 0\}$. Without loss of generality, we choose $c_{\alpha}=\frac{\alpha}{\Gamma(1-\alpha)} $ in this paper.

The next result is from \cite{PitmanYor97}.

\begin{Prop}{\rm (Pitman and Yor).}\label{p1}
Let $\{\sigma_s: s\geq 0\}$ and $\{\ga_s: s \geq 0\}$ be two
independent subordinators with respective  L\'evy measures $\alpha C
x^{-(\alpha+1)}e^{-x}dx$ and $x^{-1}e^{-x}dx$ for some $C>0$. Let
\[
\zeta(\alpha,\theta)=\frac{\ga_{\theta/\alpha}}{C\Gamma(1-\alpha)}.
\]
Then  $T=T(\alpha,\theta)=\sigma_{\zeta(\alpha,\theta)}$, and
\[
\left(\frac{V_1(T)}{T},\frac{V_2(T)}{T},\ldots\right)
\]
are independent with respective laws the $Gamma(\theta,1)$
distribution and  the two-parameter Poisson-Dirichlet distribution
$\Pi_{\alpha,\theta}$.
\end{Prop}

Let $E_{\alpha,\theta}$ denote the expectation
with respect to $\Pi_{\alpha,\theta}$. For $n \geq 1$, set
\beqn
C_{\alpha,\theta}&=&\frac{\Gamma(\theta+1)}{\Gamma(\frac{\theta}{\alpha}+1)},\label{pre0}\\
 C_{\alpha,\theta,n}
&=&
\frac{\Gamma(\theta+1)\Gamma(\frac{\theta}{\alpha}+n)\alpha^{n-1}}{\Gamma(\theta+n\alpha)
\Gamma(\frac{\theta}{\alpha}+1)\Gamma(1-\alpha)^n}.\label{pre3}
\eeqn

The following change of measure formula is obtained in \cite{PPY92}.

\begin{Prop}\label{p2}{\rm (Perman, Pitman and Yor).}
 For any bounded measurable function $f$ on ${\mathbb R}_{+}^{\infty}$,

\be\label{changefor}
E_{\alpha,\theta}(f(P_1,P_2,\ldots))=C_{\alpha,\theta}
E\left(\tau_1^{-\theta}f\left(\frac{V_1(\tau_1)}{\tau_1},
\frac{V_2(\tau_1)}{\tau_1},\ldots\right)\right),
\ee where the law of
\[
\left(\frac{V_1(\tau_1)}{\tau_1},\frac{V_2(\tau_1)}{\tau_1},\cdots\right)
\]
is $\Pi_{\alpha,0}$.
\end{Prop}

Now we are ready to derive the following distributional results.

\begin{theorem}\label{pre-t1}
For each $\beta>0$, define \be g_{\alpha,\beta}(x)=
P\left(P_1(\alpha,\beta)\leq x \right).\label{pre5} \ee Then for any
$n \geq 1$, the joint density function of
$(P_1(\alpha,\theta),\cdots,P_n(\alpha,\theta))$ is given by
\be\label{pre6}
h_{\alpha,\theta,n}(p_1,\cdots,p_n)=C_{\alpha,\theta,n}\frac{\left(1-\sum_{i=1}^n
p_i\right)^{\theta+n\alpha-1}}{(\prod_{i=1}^n p_i)^{1+\alpha}}
g_{\alpha,\theta+n\alpha}\left(\frac{p_n}{1-\sum_{i=1}^n
p_i}\right). \ee
\end{theorem}
\proof By Proposition~\ref{p2} and Perman's formula (cf.
\cite{Perman93}), for any non-negative product measurable function
$f$ and any any $n > 1$, the joint density function of
$\left(\tau_1,\frac{V_1(\tau_1)}{\tau_1},\cdots,\frac{V_n(\tau_1)}{\tau_1}\right)$
is given by \be\label{pre9}
\phi_n\left(t,p_1,\cdots,p_n\right)=(c_{\alpha})^{n-1}\hat{p}_n^{-1}\left(p_1\cdots
p_{n-1}\right)^{-(1+\alpha)}t^{-(\theta
+(n-1)\alpha)}\phi_1\left(t\hat{p}_n, p_n/\hat{p}_n\right). \ee
where \be \hat{p}_n=1-p_1-\cdots-p_{n-1}, \label{pre00} \ee
 and
$\phi_1(t,u)$ satisfies \be \label{pre11}
\phi_1(t,u)=c_{\alpha}t^{-\alpha}u^{-(1+\alpha)}\int_0^{\frac{u}{1-u}\wedge
1}\phi_1(t(1-u), v)dv. \ee

Integrating out the $t$ coordinate, it follows from \rf{pre11} that
\beqn &&h_{\alpha,\theta,n}\left(p_1,\cdots,p_n\right)\label{pre10}\\
&=&C_{\alpha,\theta}(c_{\alpha})^{n-1}\hat{p}_n^{-1}\left(p_1\cdots
p_{n-1}\right)^{-(1+\alpha)}\int_0^{\infty}t^{-(\theta
+(n-1)\alpha)}
\phi_1\left(t\hat{p}_n, p_n/\hat{p}_n\right)dt\nn\\
&=& C_{\alpha,\theta}(c_{\alpha})^{n-1}\hat{p}_n^{\theta
+(n-1)\alpha-2} \left(p_1\cdots
p_{n-1}\right)^{-(1+\alpha)}\int_0^{\infty}
s^{-(\theta +(n-1)\alpha)}\phi_1\left(s, p_n/\hat{p}_n\right)ds\nn\\
&=& C_{\alpha,\theta}(c_{\alpha})^{n}\frac{\hat{p}_n^{\theta
+n\alpha-1}}{\left(p_1\cdots p_{n-1}p_{n}\right)^{(1+\alpha)}}
\int_0^{\frac{p_n}{\hat{p}_{n+1}}\wedge 1} dx\int_0^{\infty}
s^{-(\theta +n\alpha)}\phi_1\left(s\left(1-p_n/\hat{p}_n\right), x\right)ds\nn\\
&=& C_{\alpha,\theta}(c_{\alpha})^{n}\frac{(\hat{p}_{n+1})^{\theta +n\alpha-1}}
{\left(p_1\cdots p_{n-1}p_n\right)^{(1+\alpha)}}\int_0^{\frac{p_n}{\hat{p}_{n+1}}\wedge 1} dx
\int_0^{\infty}u^{-(\theta +n\alpha)}\phi_1(u, x)du\nn\\
&=&
\frac{C_{\alpha,\theta}(c_{\alpha})^{n}}{C_{\alpha,\theta+n\alpha}}
\frac{(\hat{p}_{n+1})^{\theta +n\alpha-1}}{\left(p_1\cdots
p_{n-1}p_n\right)^{(1+\alpha)}}g_{\alpha,\theta+n\alpha}\left(\frac{p_n}{1-\sum_{i=1}^n
p_i}\right),\nn \eeqn which leads to \rf{pre6}.

\hfill $\Box$

{\bf Remark.} This result appears in Handa \cite{handa07} where a different proof is used.

\begin{theorem}\label{pre-t3} For any  $s >0$,

\beqn \nu_{\alpha,\theta,1}(s)&=&P(V_1(T)\leq
s)\label{pre23}\\
&=& \left(1+ c_\alpha
s^{-\alpha}\int_1^{\infty}z^{-(1+\alpha)}e^{- s
z}dz\right)^{-\theta/\alpha}. \nn
\eeqn
\end{theorem}
\proof  For each $s>0$, it follows from Proposition~\ref{p1} and the
property of the Poisson random measure that
\beq
\nu_{\alpha,\theta,1}(s) &=& E(P(V_1(T)\leq s |\zeta(\alpha,\theta)))\\
&=& E\left(\exp\left\{-\alpha C\zeta(\alpha,\theta)\int_s^{\infty}x^{-(\alpha+1)}e^{-x}dx\right\}\right)\\
&=& E\left(\exp\left\{-c_{\alpha}\ga_{\theta/\alpha}s^{-\alpha}
\int_1^{\infty}z^{-(\alpha+1)}e^{-sz}dz\right\}\right) \eeq which
leads to \rf{pre23}.

\hfill $\Box$

\section{Moderate Deviations for the two-parameter Poisson-Dirichlet Distribution}

By theorem 6.1 in \cite{handa07}, when $\theta$ goes to infinity
${\bf
P}(\alpha,\theta)=(P_1(\alpha,\theta),P_2(\alpha,\theta),\cdots)$
approaches a non-trivial random sequence when scaled by a factor of
$\theta$ and shifted by
$$
\beta(\alpha,\theta) =\log \theta-(\alpha+1)\log\log \theta -\log
\Gamma(1-\alpha).
$$

In \cite{Feng07}, the LDP has been established associated with the
limit
\[
\lim_{\theta \ra \infty}{\bf P}(\alpha,\theta)=(0,0,\ldots).
\]

Replacing the scaling factor by $a(\theta)$ satisfying

\be\label{mdp-seq}
\lim_{\theta\to\infty}\frac{a(\theta)}{\theta}=0,~~\lim_{\theta\to\infty}
{a(\theta)} =\infty,\ee we still have

\be\label{pd-mdp1}
\lim_{\theta \ra \infty}a(\theta)\left({\bf
P}(\alpha,\theta)-\frac{\beta(\alpha,\theta)}{\theta}(1,1,\ldots)\right)\ra
(0,0,\ldots). \ee

 The LDP associated with \rf{pd-mdp1} is called the MDP
for ${\bf
P}(\alpha,\theta)=(P_1(\alpha,\theta),P_2(\alpha,\theta)\cdots,)$.
This MDP will be established in this section through a series of
lemmas.

The first lemma establishes the MDP for $V_1(T)/\theta$.

\begin{lemma}\label{t-pd-mdp1}
The MDP holds for $V_1(T)/\theta$ with speed
$\frac{a(\theta)}{\theta}$ and rate function
\[
J_1(x)=\left\{\begin{array}{ll}
\displaystyle x,&  \displaystyle x\geq 0\\
\displaystyle\infty,& \mb{otherwise}.
\end{array}\right.
\]
\end{lemma}

\proof  For any fixed $x$, we have

\be\label{pd-mdp4}
P\left\{a(\theta)\left(\frac{V_1(T)-\beta(\alpha,\theta)}{\theta}\right)
\leq x\right\} =P \left(V_1(T)<\frac{\theta}{a(\theta)}x +
\beta(\alpha,\theta)\right) .
 \ee

Assume that $$\lim_{\theta \ra \infty}[\frac{\theta}{a(\theta)}x +
\beta(\alpha,\theta))]=+\infty.$$

Then it follows from \rf{pre23} that
$$
P \left(V_1(T)<\frac{\theta}{a(\theta)}x +
\beta(\alpha,\theta)\right) \sim
\left(1+\frac{c_\alpha}{(\frac{\theta}{a(\theta)}x +
\beta(\alpha,\theta))^{\alpha+1}}e^{-(\frac{\theta}{a(\theta)}x +
\beta(\alpha,\theta))}\right)^{-\theta/\alpha}
$$

Therefore

\beqn
&&\limsup_{\theta \ra \infty}\frac{a(\theta)}{\theta} \log
P\left(a(\theta)\left(\frac{V_1(T)-\beta(\alpha,\theta)}{\theta}\right)\leq
x\right)\nn\\
&=&\lim_{\theta \ra \infty}\frac{a(\theta)}{\theta}
\log\left(1+\frac{c_\alpha (\log
\theta)^{\alpha+1}\Gamma(1-\alpha)}{\theta(\frac{\theta}{a(\theta)}x
+\beta(\alpha,\theta))^{\alpha+1}}e^{-\frac{\theta}{a(\theta)}x}\right)^{-\theta/\alpha}
\label{appro}\\
&=&\left\{\begin{array}{ll}
      0,& x \geq 0\\
      -\infty,& x <0
      \end{array}\right..\nn
\eeqn

If there exists a subsequence $\theta'$ such that the $\lim_{\theta'
\ra \infty}(\frac{\theta'}{a(\theta')}x + \beta(\alpha,\theta'))$
exists in $[-\infty, +\infty)$, then $x$ must be strictly negative.
Since, by Theorem~\ref{pre-t3}, $V_1(T)$ converges to infinity as $\theta$ converges to infinity, it follows that

\be\label{mdp-pp2add} \limsup_{\theta' \ra
\infty}\frac{a(\theta')}{\theta'} \log
P\left(a(\theta')\left(\frac{V_1(T)-\beta(\alpha,\theta')}{\theta'}\right)\leq
x\right)=-\infty. \ee

Putting \rf{appro} and \rf{mdp-pp2add} together, one gets

\be\label{pd-mdp6} \lim_{\theta \ra \infty}\frac{a(\theta)}{\theta}
\log
P\left(a(\theta)\left(\frac{V_1(T)-\beta(\alpha,\theta)}{\theta}\right)\leq
x\right)= 0,\ x\geq 0, \ee and \be\label{mdp-pp2} \limsup_{\theta
\ra \infty}\frac{a(\theta)}{\theta} \log
P\left(a(\theta)\left(\frac{V_1(T)-\beta(\alpha,\theta)}{\theta}\right)\leq
x\right)=-\infty,\  x<0. \ee

For $x\geq 0$, it follows from \rf{pd-mdp4} and \rf{pre23} that

\beqn &&\limsup_{\theta \ra \infty}\frac{a(\theta)}{\theta} \log
P\left(a(\theta)\left(\frac{V_1(T)-\beta(\alpha,\theta)}{\theta}\right)
\geq x\right)\nn\\
&& \ \ \ \ =\limsup_{\theta \ra \infty}\frac{a(\theta)}{\theta}
\log\left(1-\left(1+\frac{c_\alpha e^{-(\frac{\theta}{a(\theta)}x +
\beta(\alpha,\theta))}}{\theta(\frac{\theta}{a(\theta)}x
+ \beta(\alpha,\theta))^{\alpha+1}}\right)^{-\theta/\alpha}\right)\label{mdp-pp1}\\
&&\ \ \ \ = \limsup_{\theta \ra \infty}\frac{a(\theta)}{\theta}
\log\left(1+\frac{c_\alpha e^{-(\frac{\theta}{a(\theta)}x +
\beta(\alpha,\theta))}}{\theta(\frac{\theta}{a(\theta)}x +
\beta(\alpha,\theta))^{\alpha+1}}\right)^{-\theta/\alpha}= -x.\nn
\eeqn

A combination of \rf{mdp-pp1} and \rf{mdp-pp2} implies that the laws
of
$a(\theta)\left(\frac{V_1(T)-\beta(\alpha,\theta)}{\theta}\right)$
is exponentially tight.

Similarly, we can get that for $x>0$ and $\delta>0$ with
$x-\delta>0$,

\beqn
&& \lim_{\theta \ra \infty}\frac{a(\theta)}{\theta}
\log P\left(a(\theta) \left(\frac{V_1(T)-\beta(\alpha,\theta)}{\theta}\right)\in (x-\delta,x+\delta)\right)\nn\\
&& \ \ \ \ = \lim_{\theta \ra \infty}\frac{a(\theta)}{\theta} \log
P\left(a(\theta)
\left(\frac{V_1(T)-\beta(\alpha,\theta)}{\theta}\right) \in
(x-\delta,x+\delta)\right)\label{fengloc1}\\
&&\ \ \ \ =-x+\delta.\nn \eeqn

The equality \rf{pd-mdp6} combined with \rf{mdp-pp2} implies that

\beqn
&&\lim_{\delta \ra 0}\limsup_{\theta \ra
\infty}\frac{a(\theta)}{\theta}
\log P\left(a(\theta) \left(\frac{V_1(T)-\beta(\alpha,\theta)}{\theta}\right)\in (-\delta,\delta)\right)\nn\\
&&\ \ \ \ =\lim_{\delta \ra 0}\liminf_{\theta \ra
\infty}\frac{a(\theta)}{\theta} \log P\left(a(\theta)
\left(\frac{V_1(T)-\beta(\alpha,\theta)}{\theta}\right) \in
(-\delta,\delta)\right)\label{fengloc2}\\
&& \ \ \ \ =0.\nn \eeqn

The lemma now follows from \rf{fengloc1}, \rf{fengloc2}, and the exponential
tightness.

\hfill $\Box$

Set
\[
\ga(\theta)=\frac{a(\theta)\beta(\alpha,\theta)}{\theta},
\]
and, without loss of generality, we can assume that
\[
\lim_{\theta \to\infty}\ga(\theta)=c \in [0,+\infty].
\]

It is clear that
\be\label{gamaspeed}
\frac{a(\theta)}{\ga^2(\theta)}=
\frac{\theta^2}{a(\theta)\beta^2(\alpha,\theta)}\ra \infty, \ \
\theta \ra \infty. \ee

If $c <\infty$, it follows from Corollary~3.1 in \cite{FengGao08}
that for any $L >0$

\be\label{exponential} \limsup_{\theta\to
\infty}\frac{a(\theta)}{\theta}\log
P\left\{\ga(\theta)|\frac{\theta}{T}-1|\geq
L\right\}=-\infty. \ee

For  $c=\infty $, and any $1>\delta >0$ \be\label{exponential0}
\left\{\ga(\theta)|\frac{\theta}{T}-1|\geq L\right\}
\subset \left\{\ga(\theta)|\frac{T}{\theta}-1|\geq
L(1-\delta)\right\} \bigcup
\left\{|\frac{T}{\theta}-1|\geq \delta\right\}. \ee

Since $\ga(\theta)\leq \beta(\alpha,\theta)$ for large $\theta$ and $\lim_{\theta \ra
\infty}\frac{\beta(\alpha,\theta)}{\sqrt{\theta}}=0$, it follows
from the MDP (Theorem~3.2 in \cite{FengGao08}) for
$T/\theta$, and \rf{gamaspeed} that
  \be\label{exponential1}
\begin{aligned}
&\limsup_{\theta \ra \infty}\frac{a(\theta)}{\theta}\log
P\left\{\ga(\theta)|\frac{T}{\theta}-1|\geq
(1-\delta)L\right\}\\
=&\limsup_{\theta\to \ra
\infty}\frac{a(\theta)}{\ga^2(\theta)}\frac{\ga^2(\theta)}{\theta}\log
P\left\{\ga(\theta)|\frac{T}{\theta}-1|\geq
(1-\delta)L\right\}=-\infty,
\end{aligned}
\ee which combined with Corollary~3.1  in \cite{FengGao08} and
\rf{exponential0} shows that \rf{exponential} still holds in this
case. Therefore
$a(\theta)(P_1(\alpha,\theta)-\frac{\beta(\alpha,\theta)}{\theta})$
and $\frac{\theta}{T}a(\theta)
(\frac{V_1(T)-\beta(\alpha,\theta)}{\theta})$ are
exponentially equivalent.

Since
$\frac{\theta}{T}a(\theta)(\frac{V_1(T)
-\beta(\alpha,\theta)}{\theta})$ is exponentially equivalent to
$a(\theta)(\frac{V_1(T)-\beta(\alpha,\theta)}{\theta})$
by Lemma~2.1 and Corollary~3.1 in \cite{FengGao08}, it follows that
$a(\theta)(P_1(\alpha,\theta)-\frac{\beta(\alpha,\theta)}{\theta})$
and
$a(\theta)(\frac{V_1(T)-\beta(\alpha,\theta)}{\theta})$
are exponentially equivalent. Thus we have the following result.

\begin{lemma}\label{P-pd-mdp1}
The MDP holds for $P_1(\alpha,\theta)$ with speed
$\frac{a(\theta)}{\theta}$ and rate function
\[
J_1(x)=\left\{\begin{array}{ll}
\displaystyle x,&  \displaystyle x\geq 0\\
\displaystyle\infty,& \mb{otherwise}.
\end{array}\right.
\]
\end{lemma}

For each $n \geq 2$, we have

\begin{lemma}\label{mdp-Pk}
 The family $ \Big\{P\Big(a(\theta)\Big(P_1(\alpha,\theta)-\frac{\beta(\alpha,\theta)}{\theta}\cdots,
  P_n(\alpha,\theta)-\frac{\beta(\alpha,\theta)}{\theta}
\Big)\in\cdot \Big):~ \theta>0 \Big\} $ satisfies a LDP on ${\mathbb
R^n}$ with speed $\frac{a(\theta)}{\theta} $ and rate function
\be\label{rate-PDk-I} I_n(x_1,\cdots,x_n) = \left\{\begin{array}{ll}
\displaystyle \sum_{i=1}^n x_i,&  if\quad 0\leq x_n\leq \cdots\leq x_1.\\
\displaystyle +\infty,& otherwise.
\end{array}\right.
\ee
\end{lemma}

\proof \quad  It follows from \rf{pd-mdp4} that for $x_1\geq
x_2\cdots\geq x_n$  and
$\frac{\theta}{a(\theta)}x_n+\beta(\alpha,\theta)>0$, the density
function $g_{\alpha,\theta,n}(x_1,\ldots,x_n)$ of
$a(\theta)\Big(P_1(\alpha,\theta)-\frac{\beta(\alpha,\theta)}{\theta}\cdots,
  P_n(\alpha,\theta)-\frac{\beta(\alpha,\theta)}{\theta}
\Big)$ is
\beqn
&&g_{\alpha,\theta,n}(x_1,\ldots,x_n)=
\left(\frac{1}{a(\theta)}\right)^nC_{n,\alpha,\theta}\left(\prod_{i=1}^n
\left(\frac{\theta}{\frac{\theta}{a(\theta)}x_i
+\beta(\alpha,\theta)}\right)^{\alpha+1}\right) \label{multi-mdp1}\\
&&\ \ \times\left(1-\left(\frac{\theta}{a(\theta)}\sum_{i=1}^n x_i +n
 \beta(\alpha,\theta)\right)/\theta\right) ^{\theta+n \alpha -1}
g_{\alpha,\theta+n \alpha}\left(\frac{\frac{\theta}{a(\theta)} x_n +
 \beta(\alpha,\theta)}{\theta-\left(
\frac{\theta}{a(\theta)}\sum_{i=1}^n x_i +n
 \beta(\alpha,\theta)\right)}\right).\nn
\eeqn
By Theorem~\ref{pre-t3} and direct calculation, for $x_n>0$
$$
\frac{a(\theta)}{\theta}\log g_{\alpha,\theta+n
\alpha}\left(\frac{\frac{\theta}{a(\theta)} x_n +
 \beta(\alpha,\theta)}{\theta-\left(
\frac{\theta}{a(\theta)}\sum_{i=1}^n x_i +n
 \beta(\alpha,\theta)\right)}\right)\to 0.
$$

For $x_n<0$, set
$$
\psi(n,x,\theta,\alpha)=a(\theta)\left(\frac{\frac{\theta}{a(\theta)}
x_n +
 \beta(\alpha,\theta)}{\theta-\left(
\frac{\theta}{a(\theta)}\sum_{i=1}^n x_i +n
 \beta(\alpha,\theta)\right)}
 -\frac{\beta(\alpha,\theta+n\alpha)}{\theta+n\alpha}\right) .
$$
Then
$$
\begin{aligned}
& g_{\alpha,\theta+n \alpha}\left(\frac{\frac{\theta}{a(\theta)} x_n
+
 \beta(\alpha,\theta)}{\theta-\left(
\frac{\theta}{a(\theta)}\sum_{i=1}^n x_i +n
 \beta(\alpha,\theta)\right)}\right)\\
=& P\left(a(\theta)\left(P_1(\alpha,\theta+n\alpha)
-\frac{\beta(\alpha,\theta+n\alpha)}{\theta+n\alpha}\right)
<\psi(n,x,\theta,\alpha)\right)
\end{aligned}
$$
and
$$
\lim_{\theta\to\infty}\psi(n,x,\theta,\alpha)=x_n<0
$$
which implies that
$$
\lim_{\theta\to\infty}\frac{a(\theta)}{\theta}\log
g_{\alpha,\theta+n \alpha}\left(\frac{\frac{\theta}{a(\theta)} x_n +
 \beta(\alpha,\theta)}{\theta-\left(
\frac{\theta}{a(\theta)}\sum_{i=1}^n x_i +n
 \beta(\alpha,\theta)\right)}\right)= -\infty.
$$

Therefore
\beqn
\frac{a(\theta)}{\theta}\log g_{\alpha,\theta,n}(x_1,\ldots,x_n)\ra -\sum_{i=1}^n x_i, \ x_n >0, \label{multi-mdp2}\\
\frac{a(\theta)}{\theta}\log g_{\alpha,\theta,n}(x_1,\ldots,x_n)\ra
-\infty,\  x_n <0, \label{multi-mdp3}. \eeqn

For $ x_1\geq x_2\cdots\geq x_n$, let $B((x_1,\ldots,x_n), \delta)$
denote the closed ball centered at $(x_1,\ldots,x_n)$ with radius
$\delta$, and $B^{\circ}((x_1,\ldots,x_n), \delta)$ be the
corresponding open ball. Then for $x_n >0$,

\beqn &&\lim_{\delta \ra 0}\limsup_{\theta \ra
\infty}\frac{a(\theta)}{\theta} \log
P\left(a(\theta)\Big(P_1(\alpha,\theta)-\frac{\beta(\alpha,\theta)}{\theta}\cdots,
  P_n(\alpha,\theta)-\frac{\beta(\alpha,\theta)}{\theta}
\Big)\in B((x_1,\ldots,x_n), \delta)\right)\nn\\
&&  =\lim_{\delta \ra 0}\liminf_{\theta \ra
\infty}\frac{a(\theta)}{\theta} \log
P\left(a(\theta)\Big(P_1(\alpha,\theta)-\frac{\beta(\alpha,\theta)}{\theta}\cdots,
  P_n(\alpha,\theta)-\frac{\beta(\alpha,\theta)}{\theta}
\Big)\in B^{\circ}((x_1,\ldots,x_n), \delta)\right)\nn\\
&&  =-\sum_{i=1}^n x_i,\nn \eeqn

and for any $x_n<0$,

\beqn &&\lim_{\delta \ra 0}\limsup_{\theta \ra
\infty}\frac{a(\theta)}{\theta} \log
P\left(a(\theta)\Big(P_1(\alpha,\theta)-\frac{\beta(\alpha,\theta)}{\theta}\cdots,
  P_n(\alpha,\theta)-\frac{\beta(\alpha,\theta)}{\theta}
\Big)\in B((x_1,\ldots,x_n), \delta)\right)\nn\\
&& =\lim_{\delta \ra 0}\liminf_{\theta \ra
\infty}\frac{a(\theta)}{\theta} \log
P\left(a(\theta)\Big(P_1(\alpha,\theta)-\frac{\beta(\alpha,\theta)}{\theta}\cdots,
  P_n(\alpha,\theta)-\frac{\beta(\alpha,\theta)}{\theta}
\Big)\in B^{\circ}((x_1,\ldots,x_n), \delta)\right)\nn\\
&&  =-\infty,\nn \eeqn

If $x_{r-1}>0, x_r=0$ for some $1 \leq r \leq n
$, then the upper estimate is obtained from that of
$a(\theta)\Big(P_1(\alpha,\theta)-\frac{\beta(\alpha,\theta)}{\theta}\cdots,
  P_{r-1}(\alpha,\theta)-\frac{\beta(\alpha,\theta)}{\theta}
\Big)$.  The lower estimates when $x_r=0$ for some $1 \leq r \leq n
$ are obtained by approximating the boundary with open subsets away from the boundary.

Noting that $\bigcup_{i=1}^n
\{a(\theta)\Big(P_i(\alpha,\theta)-\frac{\beta(\alpha,\theta)}{\theta}\Big)
>L\} =\{a(\theta)\Big(P_1(\alpha,\theta)-\frac{\beta(\alpha,\theta)}{\theta}\Big) >L\}$,
it follows that

\be\label{appro2} \lim_{L\ra \infty}\limsup_{\theta \ra
\infty}\frac{a(\theta)}{\theta}\log P\left\{\bigcup_{i=1}^n
\left\{a(\theta)\Big(P_i(\alpha,\theta)-\frac{\beta(\alpha,\theta)}{\theta}\Big)
>L\right\}\right\}=-\infty. \ee

On the other hand,

\beqn &&\limsup_{\theta \ra \infty}\frac{a(\theta)}{\theta}\log
P\left\{\bigcup_{i=1}^n\left\{a(\theta)\Big(P_i(\alpha,\theta)-\frac{\beta(\alpha,\theta)}{\theta}\Big) <-L\right\}\right\}\label{appro3}\\
&&\ \ \ \ \leq \limsup_{\theta \ra
\infty}\frac{a(\theta)}{\theta}\log
P\left\{a(\theta)\Big(P_i(\alpha,\theta)-\frac{\beta(\alpha,\theta)}{\theta}\Big)
\leq -L\right\}=-\infty. \nn\eeqn

These lead to the exponential tightness and the lemma.

\hfill $\Box$

Now we are ready to establish the MDP for $ (P_1(\alpha,\theta),P_2(\alpha,\theta),\ldots)$.

\begin{theorem}\label{mdp-PD}
 For each $n \geq 1$, the family
 $ \Big\{P\Big(a(\theta)\Big( P_1(\theta)-\frac{\beta(\theta)}{\theta},\cdots, P_n(\theta)
 -\frac{\beta(\theta)}{\theta},\cdots
\Big)\in\cdot \Big):~ \theta>0 \Big\} $ satisfies a LDP on ${\mathbb
R^\infty}$ with speed $\frac{a(\theta)}{\theta} $ and rate function
\be\label{rate-PD} I(x_1,x_2,\cdots) = \left\{\begin{array}{ll}
\displaystyle\sum_{i=1}^{\infty}x_i,&  \displaystyle x_1\geq\cdots\geq 0 \\
\displaystyle\infty,& \mb{otherwise}.
\end{array}\right.
\ee
\end{theorem}
\proof Identify ${\mathbb R^\infty}$ with the projective limit of
${\mathbb R^n}, n=1,\ldots.$ Then the theorem follows from Theorem
3.3 in \cite{DaGa87} and Lemma~\ref{mdp-Pk}.

\hfill $\Box$

\section{Moderate Deviations for the Homozygosity}

For each $m \geq 2$, it was shown in \cite{handa07} that the scaled
homozygosity
$$\sqrt{\theta}[\frac{\theta^{m-1}\Gamma(1-\alpha)}{\Gamma(m-\alpha)}H_m({\bf
P}(\alpha,\theta))-1] \Rightarrow Z_{\alpha,m}$$ where $Z_{\alpha,m}$ is a normal random variable with mean zero and
variance
\[
\sigma^2_{\alpha,m}=\frac{\Gamma(2m-\alpha)\Gamma(1-\alpha)}{\Gamma(m-\alpha)^2}+\alpha-m^2.
\]

It is thus natural to
consider the MDP for $\frac{\theta^{m-1}}{\Gamma(m)}H_m({\bf
P}(\alpha,\theta))$ or equivalently the LDP for the family
$\{a(\theta)[\frac{\theta^{m-1}\Gamma(1-\alpha)}{\Gamma(m-\alpha)}H_m({\bf
P}(\alpha,\theta))-1]:\theta
>0\}$
 for a scale $a(\theta)$ satisfying
\be\label{mdpscale2} \lim_{\theta \ra \infty}a(\theta)=\infty,
\lim_{\theta \ra \infty}\frac{a(\theta)}{\sqrt{\theta}}=0, \ee
which is different from \rf{mdp-seq}.

The MDP in the case of $\alpha=0$ has been established in \cite{FengGao08} where that the following additional
restriction on $a(\theta)$ is used: for
some $0<\epsilon<1/(2m-1)$,

\be\label{a-speed}
\liminf_{\theta\to\infty}\frac{a^{1-\varepsilon}(\theta)}{\theta^{(m-1)/(2m-1)}}>0.
\ee

This condition is also needed for the two-parameter model. As shown in \cite{FengGao08}, the conditions
 \rf{mdpscale2} and \rf{a-speed} guarantee that there exist $\tau>0$
positive integer $l\geq 3\vee \frac{2}{(2m-1)\varepsilon}$, and $r(\theta) $  that
grows faster than a positive power of
$\theta$ such that
$$
\lim_{\theta\to\infty}\frac{a (\theta)}{\theta^{\tau}}=+\infty
$$
and
$$
\lim_{\theta\to\infty}
\frac{r(\theta)^{m-1}}{a^{(l-2)/l}(\theta)} =0,
~~~\lim_{\theta\to\infty}\frac{ a^2(\theta)r(\theta)}{\theta}
=\infty.
$$

For any $n\geq1$, set
\beq
G_{\alpha,\theta,r}^{(n)}&=&\sum_{i=1}^{\infty}V^n_i(T)I_{\{V_i(T)\leq r(\theta)\}},\\
G_{\alpha,\theta }^{(n)}&=& \sum_{i=1}^{\infty}V^n_i(T),
\eeq
and
\[
G_{\alpha,\theta,r}=\left(G_{\alpha,\theta,r
}^{(1)}-E(G_{\alpha,\theta,r }^{(1)}), G_{\alpha,\theta,r
}^{(m)}-E(G_{\alpha,\theta }^{(m)})\right).
\]

For any $s,t$ in ${\mathbb R}$, define
\[
\Lambda(s,t)=\frac{1}{2}\left(s^2+\frac{2\Gamma(m-\alpha)\Gamma(m+1)}{
\Gamma(m)\Gamma(1-\alpha)}st +(\frac{\Gamma (2m-\alpha) }{\Gamma
(1-\alpha)}+ \alpha (\frac{\Gamma(m-\alpha)}{\Gamma(1-\alpha)})^2)t^2\right).
\]

It follows by direct calculation that the Fenchel-Legendre transform of $\Lambda(s,t)$ is given by
\beqn
\Lambda^*(x,y) &=& \sup_{s,t}\{sx+ty -\Lambda(s,t)\}\nn\\
&=&\frac{\Gamma(1-\alpha)}{2(
\Gamma(1-\alpha)\Gamma(2m-\alpha) +(\alpha-m^2)
\Gamma^2(m-\alpha))}\label{trunrate}\\
&&\times\left( (\Gamma (2m-\alpha)+\alpha \frac{\Gamma^2(m-\alpha)}{\Gamma(1-\alpha)})
x^2 -2m\Gamma(m-\alpha)x y +
\Gamma(1-\alpha)y^2\right),\nn
\eeqn
for $x,y$ in ${\mathbb R}$.

\begin{lemma}\label{homozygosity-mdp-lem}
The family $\{\frac{a(\theta)}{\theta}
G_{\alpha,\theta,r}:\theta
>0\}$ satisfies a LDP on space ${\mathbb R}^2$ with speed $\frac {a^2(\theta)}{\theta}$ and
rate function $\Lambda^*(\cdot,\cdot)$.
\end{lemma}
\proof For any $s,t\in\mathbb R$, let
\[
g(x)=sx + t x^m
\]
and
\[
\varphi_{r}(x)=\frac{g(x)I_{\{x\leq r(\theta)\}}}{a(\theta)}.
\]

It follows by direct calculation that
\be\label{ranenviro1}
\begin{aligned}
&\int_0^{r(\theta)}(e^{\varphi_{r}(x)}-1) x^{-(1+\alpha)}e^{-x}dx\\
=&\int_0^{r(\theta)}\frac{g(x)}{a(\theta)} x^{-(1+\alpha)}e^{-x}dx
 +\frac{1}{2}\int_0^{r(\theta)}\frac{g^2(x)}{a^2(\theta)}
x^{-(1+\alpha)}e^{-x}dx\\
&+\sum_{k=3}^l \frac{1}{k!}
\frac{1}{a^{k}(\theta)}\int_0^{\gamma(\theta)}|s x+t x^m
|^k x^{-(1+\alpha)} e^{-x} dx\\
&+O\left(\sum_{k=l+1}^\infty \frac{1}{k!} \frac{1}{a^{k}(\theta)}
(|s|+|t|\gamma(\theta) ^{m-1})^k\Gamma(k-\alpha)\right)\\
=&\int_0^{r(\theta)}\frac{g(x)}{a(\theta)} x^{-(1+\alpha)}e^{-x}dx
 +\frac{1}{2}\int_0^{r(\theta)}\frac{g^2(x)}{a^2(\theta)}
x^{-(1+\alpha)}e^{-x}dx+o\left(\frac{1}{a^2(\theta)}\right),
\end{aligned}
\ee
which implies that for $\theta$ large enough,
\[
|\int_0^{r(\theta)}(e^{\varphi_{r}(x)}-1)
x^{-(1+\alpha)}e^{-x}dx| <c^{-1}_{\alpha}.
\]

By the Campbell's theorem we get that
\beqn
&&E\left(\exp\left\{\frac{1}{a(\theta)}(sG_{\alpha,\theta,r}^{(1)}
+tG_{\alpha,\theta,r}^{(m)})\right\}\right)\nn\\
&&\ \ \ \ \ =E\left(\exp\left\{\sum_{i=1}^\infty
\varphi_{r}(V_i(T))\right\}\right)\nn\\
&&\ \ \ \ \ = E\left(E\left( \exp\left\{\sum_{i=1}^\infty
\varphi_{r}(V_i(T))\right\}|\zeta(\alpha,\theta)\right) \right)\label{ranenviro2}\\
&&\ \ \ \ \ =E \left(\exp\{
c_{\alpha}\ga(\frac{\theta}{\alpha})\int_0^{r(\theta)}(e^{\varphi_{r}(x)}-1)
x^{-(1+\alpha)}e^{-x}dx\}\right)\nn\\
&&\ \ \ \ \ =\exp\left\{-\frac{\theta}{\alpha}\log\left(1- c_{\alpha}\int_0^{r(\theta)}
(e^{\varphi_{r}(x)}-1)
x^{-(1+\alpha)}e^{-x}dx\right)\right\}.\nn
\eeqn

Putting \rf{ranenviro1} and \rf{ranenviro2} together, we get that

\beq &&
E\left(\exp\left\{\frac{1}{a(\theta)}(s(G_{\alpha,\theta,r}^{(1)}-E(G_{\alpha,\theta,r}^{(1)}))
+t(G_{\alpha,\theta,r}^{(m)}-E(G_{\alpha,\theta,r}^{(m)})))\right\}\right)\nn\\
&&\ \ =\exp\left\{\frac{\theta c_{\alpha}}{2\alpha
a^2(\theta)}\left(c_{\alpha}(\int_0^{\infty}g(x)x^{-(1+\alpha)}
e^{-x}dx)^2  + \int_{0}^{\infty}g^2(x)x^{-(1+\alpha)}
e^{-x}dx+ o(\frac{1}{a^2(\theta)})\right)\right\}\nn\\
&&\ \  =  \exp\left(\frac{\theta}{a^2(\theta)}(\Lambda(s,t)+
o(\frac{1}{a^2(\theta)}))\right),\nn \eeq which leads to
\begin{equation}
\lim_{\theta \ra \infty}\frac{a^2(\theta)}{\theta}\log
E\left(\exp\left\{\frac{1}{a(\theta)}[s(G_{\alpha,\theta,r}^{(1)}-E(G_{\alpha,\theta,r}^{(1)}))
+t(G_{\alpha,\theta,r}^{(m)}-E(G_{\alpha,\theta,r}^{(m)}))]\right\}\right)=\Lambda(s,t).
\end{equation}

The lemma now follows from \rf{trunrate} and
 the G\"artner-Ellis theorem.

\hfill $\Box$

\begin{lemma}\label{homozygosity-mdp-lem-2}
Set
$$
 G_{\alpha,\theta}=\left(T-\theta, G_{\alpha,\theta
}^{(m)}-E(G_{\alpha,\theta}^{(m)})\right).
$$
Then the family $\{\frac{a(\theta)}{\theta} G_{\alpha,\theta}:\theta
>0\}$ satisfies a LDP with speed $\frac {a^2(\theta)}{\theta}$ and
the rate function $\Lambda^*(x,y)$ .

\end{lemma}
\proof By definition for any $n\geq 1$ and any $\delta>0$,

$$
\begin{aligned}
&\limsup_{\theta\to\infty}\frac{a^2(\theta)}{\theta}\log
P\left(\left|
G_{\alpha,\theta,r}^{(m)}- G_{\alpha,\theta}^{(m)} \right|
\geq \delta \frac{\theta} {a(\theta)}\right)\\
\leq&\limsup_{\theta\to\infty}\frac{a^2(\theta)}{\theta}\log P\left(
V_1(T)\geq r(\theta)\right)\\
=&\limsup_{\theta\to\infty}\frac{a^2(\theta)}{\theta}
\log \left(1-\left(1+\frac{c_\alpha}{r^\alpha(\theta)}\int_1^\infty
z^{-(1+\alpha)}e^{-r(\theta)z}dz\right)^{-\theta/\alpha}\right)\\
\leq &\limsup_{\theta\to\infty}\frac{a^2(\theta)}{\theta} \log
\left(1-\left(1+\frac{c_\alpha}{r^{(1+\alpha)}(\theta)e^{r(\theta)}}
\right)^{-\theta/\alpha}\right)\\
=&\limsup_{\theta\to\infty}\frac{a^2(\theta)}{\theta}\log\left(\frac{\theta}{\alpha}
\log \left(1+\frac{c_\alpha}{r^{(1+\alpha)}(\theta)e^{r(\theta)}}
\right)\right)\\
\leq
&-\limsup_{\theta\to\infty}\frac{a^2(\theta)r(\theta)}{\theta}\left(1-\frac{\log
\theta}{r(\theta)}\right) \\
=&-\infty.
\end{aligned}
$$
which implies that $\frac{a(\theta)}{\theta}
{G}_{\alpha,\theta,r}$ and
$\frac{a(\theta)}{\theta}G_{\alpha,\theta}$ are exponentially
equivalent. Therefore $ \left(\frac{a(\theta)}{\theta}
G_{\alpha,\theta},\frac {a^2(\theta)}{\theta}, \Lambda^*\right) $
satisfies LDP.

\hfill $\Box$

Now we are ready to prove the main result of this section.

\begin{theorem}\label{homozygosity-mdp-thm}
The family $
a(\theta)\left(\frac{\theta^{m-1}\Gamma(1-\alpha)}{\Gamma(m-\alpha)}H_m({\bf
P}(\alpha,\theta))-1\right)$ satisfies a LDP with speed $ \frac
{a^2(\theta)}{\theta}$ and rate function $ \frac{z^2}{2\sigma^2_{\alpha,m}}$.
\end{theorem}

\proof  By direct calculation, \beq
&&a(\theta)\left(\frac{\theta^{m-1}\Gamma(1-\alpha)}{\Gamma(m-\alpha)}
H_m({\bf P}(\alpha,\theta))-1\right)\\
&&\hspace{1cm}= {a(\theta)}\left(\frac{\theta ^{m-1}
G^{(m)}_{\alpha,\theta}}{T^m
\Gamma(m-\alpha)/\Gamma(1-\alpha)}-1\right)\\
&&\hspace{1cm}= {a(\theta)} \left(\left(\frac{\theta}{T}\right)^m-1
\right)+\left(\frac{\theta}{T}\right)^m
\frac{a(\theta)(G_{\alpha,\theta}^{(m)}-E(G_{\alpha,\theta}^{(m)})}
{\theta\Gamma(m-\alpha)/\Gamma(1-\alpha)}\\
&&\hspace{1cm}=\frac{a(\theta)}{\theta} (\theta-T)
\sum_{k=1}^m\left(\frac{\theta}{T}\right)^k+\left(\frac{\theta}{T}\right)^m
\frac{a(\theta)(G_{\alpha,\theta}^{(m)}-E(G_{\alpha,\theta}^{(m)})}
{\theta\Gamma(m-\alpha)/\Gamma(1-\alpha)}.
\eeq

Noting that for any $i\geq 1$ and for any $\delta>0$,
$$
\lim_{\theta\to\infty}\frac{a^2(\theta)}{\theta}\log
P\left(\left|\left(\frac{\theta}{T}\right)^i-1\right|\geq
\delta\right)=-\infty.
$$

It then follows that
$$a(\theta)\left(\frac{\theta^{m-1}\Gamma(1-\alpha)}{\Gamma(m-\alpha)}
H_m({\bf P}(\alpha,\theta))-1\right)$$
and
$$\frac{a(\theta) m(\theta-T)}{\theta}+
\frac{a(\theta)(G_{\alpha,\theta}^{(m)}-E(G_{\alpha,\theta}^{(m)}))}
{\theta\Gamma(m-\alpha)/\Gamma(1-\alpha)}$$ are exponentially
equivalent, and so they have the same LDP.

The fact that
$$
\inf_{\frac{y\Gamma(1-\alpha)}{\Gamma(m-\alpha)}-mx=z}
\Lambda^*(x,y)= \frac{z^2}{2\sigma^2_{\alpha,m}},
$$
combined with Lemma \ref{homozygosity-mdp-lem-2} and the contraction
principle implies the theorem.

\hfill $\Box$

\section{LDP for Small Parameters}

Let
\[
\nabla =\left\{{\bf p}=(p_1,p_2,\ldots): p_1\geq p_2\geq \cdots\geq
0, \sum_{i=1}^{\infty}p_i \leq 1\right\}
\]
be equipped with the subspace topology of $[0,1]^{\infty}$,
and $M_1(\nabla)$ be the space of all probability measures on $\nabla$ equipped with
the weak topology. Then $\Pi_{\alpha,\theta}$ belongs to  $M_1(\nabla)$.

For any $\delta >0$, it follows from the GEM representation \rf{GEM1} that
\[
P\left(X_1^{\alpha,\theta}>1-\delta\right) \leq
P\left(P_1(\alpha,\theta)>1-\delta\right).
\]

By direct calculation, we have
\[
\lim_{\alpha+\theta \ra
0}P\left(X_1^{\alpha,\theta}>1-\delta\right)=1.
\]
Therefore, $\Pi_{\alpha,\theta}$ converges in $M_1(\nabla)$ to
$\delta_{(1,0,\ldots)}$ as $\alpha+\theta$
converges to zero. In this section, we establish the LDP
associated with this limit. This is a two-parameter
generalization to the result
in \cite{Feng08}.

For any $n \geq 1$, set
      \beq
      &&\nabla_n =\left\{(p_1,...,p_n,0,0,...) \in \nabla: \sum_{i=1}^n p_i = 1\right\},\\
      &&\nabla_{\infty}= \bigcup_{i=1}^{\infty}\nabla_i,
      \eeq
and

\[
a(\alpha,\theta)= \alpha \vee |\theta|,\  b(\alpha,\theta)= (-\log(a(\alpha,\theta))^{-1}.
\]

Then we have
\begin{Lem}\label{small-P}
      The family of laws of $\{P_1(\alpha,\theta):\alpha+\theta >0, 0<\alpha<1\}$ satisfies a LDP
      on $[0, 1]$ as $a(\alpha,\theta)$ goes to zero with speed
      $b(\alpha,\theta)$ and rate function
       \be\label{small-rate1}
      S_1(p)=
      \left\{\begin{array}{ll}
      0,& p= 1\\
      k,& p\in [\frac{1}{k+1},\frac{1}{k}), k=1,2,\ldots\\
      \infty,& p=0.
      \end{array}\right.
      \ee
      \end{Lem}

\proof Let $\{X_i^{\alpha,\theta}: i=1,2,\ldots\}$ be defined in \rf{GEM1}. For any $n \geq 1$, set

\[
\tilde{P}^n_1(\alpha,\theta)=\max\{X_i^{\alpha,\theta}: 1\leq i\leq n\}.
\]
Then it follows from direct calculation that for any $\delta >0$

\beq
P\{P_1(\alpha,\theta)-\tilde{P}_1^n(\alpha,\theta)>\delta\}&\leq& P\{(1-U_1)\cdots(1-U_n)\geq \delta\}\\
&\leq & \delta^{-1}\prod_{i=1}^n \frac{\theta +i\alpha}{\theta +i\alpha +1-\alpha},
\eeq
which leads to
\[
\limsup_{a(\alpha,\theta)\ra 0}b(\alpha, \theta)\log P\{P_1(\alpha,\theta)-\tilde{P}_1^n(\alpha,\theta)>\delta\}
\leq -n.
\]

Thus the families $\{\tilde{P}^n_1(\alpha,\theta): 0<\alpha<1, \theta +\alpha >0\}_{n=1,2,\ldots}$ are exponential good approximations to the family $\{P_1(\alpha,\theta): 0<\alpha<1, \theta +\alpha >0\}$. By the contraction principle, the family
     $\{\tilde{P}^n_1(\alpha,\theta): 0<\alpha<1, \theta +\alpha >0\}$ satisfies a LDP
      on $[0, 1]$ as $a(\alpha,\theta)$ goes to zero with speed
      $b(\alpha,\theta)$ and rate function

\[
      I_n(p)= \left\{\begin{array}{ll}
      0,& p =1\\
      k,& p \in [\frac{1}{k+1},\frac{1}{k}), k=1,2,\ldots,n-1\\
      n,& \mb{else.}
      \end{array}\right.
      \]

The lemma now follows from the fact that
\[
S_1(p)=\sup_{\delta>0}\liminf_{n \ra \infty}\inf_{|q-p|<\delta}I_n(q).
\]
\hfill $\Box$

\begin{theorem}\label{P3}
      The family $\{\Pi_{\alpha,\theta}:\alpha+\theta >0, 0<\alpha<1\}$ satisfies a LDP
      on $\nabla$ as $a(\alpha,\theta)$ goes to zero with speed
      $b(\alpha,\theta)$ and rate function

      \be\label{rate function3}
      S({\bf p})=
      \left\{\begin{array}{ll}
      n-1,& {\bf p} \in \nabla_n, p_n >0, n\geq 1\\
      \infty,& {\bf p}\not\in \nabla_{\infty}.
      \end{array}\right.
      \ee
      \end{theorem}
      \proof It suffices to establish the LDP for finite dimensional marginal
      distributions since
      the infinite dimensional LDP can be derived from the finite dimensional
      LDP through approximation.
      For any $n \geq 2$,  $(P_1(\alpha,\theta), P_2(\alpha,\theta),\ldots£¬
      P_n(\alpha,\theta))$ and
      $(P_1(0,\alpha+\theta), P_2(0,\alpha+\theta),\ldots£¬P_n(0,\alpha+\theta))$
      have respective joint density functions
      \[
      h_{\alpha,\theta,n}(p_1,\cdots,p_n)=C_{\alpha,\theta,n}\frac{\left(1-\sum_{i=1}^n
p_i\right)^{\theta+n\alpha-1}}{\prod_{i=1}^n p_i} P\left(P_1(\alpha,
n\alpha+\theta)\leq \frac{p_n}{1-\sum_{i=1}^n p_i}\right),
\]
      and

      \[
       g_{\alpha+\theta,n}=(\alpha+\theta)^{n}\frac{\left(1-\sum_{i=1}^n
p_i\right)^{\theta+\alpha-1}}{\prod_{i=1}^n
p_i}P\left(P_1(0,\alpha+\theta)\leq \frac{p_n}{1-\sum_{i=1}^n
p_i}\right)
      \]

Since $\lim_{a(\alpha,\theta)\ra 0}b(\alpha,\theta)\log(\alpha+\theta)=-1$ and
$\lim_{a(\alpha,\theta)\ra 0}b(\alpha,\theta)C_{\alpha,\theta,n}=-n$, it follows from
Lemma 2.4 in \cite{Feng08} and Lemma~\ref{small-P} that the family of laws of
$(P_1(\alpha,\theta), P_2(\alpha,\theta),\ldots£¬P_n(\alpha,\theta))$
satisfies a LDP as $a(\alpha,\theta)$ goes to zero with speed
      $b(\alpha,\theta)$ and rate function

      \be\label{small-fin-rate}
     S_n(p_1,...,p_n)=
      \left\{\begin{array}{ll}
       0,& (p_1,p_2,...,p_n)=(1,0...,0)\\
       l-1,& 2\leq l \leq n, \sum_{k=1}^l p_k =1, p_l >0 \\
       n + S_1\left(\frac{p_n}{1-\sum_{i=1}^n p_i}\wedge 1\right),&  \sum_{k=1}^n p_k <1, p_n>0\\
       \infty,& \mb{else}.
      \end{array}\right.
      \ee

      \hfill $\Box$

      \section*{Acknowledgement}
     Fuqing Gao would like to thank the Department of
     Mathematics and Statistics at McMaster University for their
     hospitality during his visit.

\end{document}